\documentclass[reqno, 12pt]{article}

\pdfoutput=1

\usepackage{enumerate}
\usepackage{latexsym}
\usepackage[centertags]{amsmath}
\usepackage{amsfonts}
\usepackage{amssymb}
\usepackage{amsthm}
\usepackage{newlfont}
\usepackage{graphics}
\usepackage{color}
\usepackage{float}
\usepackage{diagbox}
\usepackage{hyperref}
\usepackage{longtable}
\usepackage{rotating}
\usepackage{multirow}
\usepackage{extarrows}
\usepackage[sort,compress,numbers]{natbib}
\usepackage[utf8]{inputenc}

\newtheorem{thm}{Theorem}[section]

\theoremstyle{mydefinition}
\newtheorem{dfn}[thm]{Definition}
\theoremstyle{myremark}

\newtheorem{prob}[thm]{Open Problem}
\allowdisplaybreaks[4]

\title{Order Polytopes of Dimension $\leq 13$ are Ehrhart Positive}

\author{Feihu Liu$^{1}$, Guoce Xin$^{2,}$\thanks{This work is partially supported by the National Natural Science Foundation of China (No.12071311).},\ and Zihao Zhang$^{3}$
\\[2mm]
{\small $^{1, 2, 3}$ School of Mathematical Sciences,}\\[-0.8ex]
{\small Capital Normal University, Beijing, 100048, P.R.~China}\\
{\small $^1$ Email address: liufeihu7476@163.com}\\
{\small $^2$ Email address: guoce\_xin@163.com}\\
{\small $^3$ Email address: zihao-zhang@foxmail.com}
}

\date{December 10, 2024}

\begin{document}

\maketitle

\begin{abstract}
The order polytopes arising from the finite poset were first introduced and studied by Stanley. For any positive integer $d\geq 14$, Liu and Tsuchiya proved that there exists a non-Ehrhart positive order polytope of dimension $d$. They also proved that any order polytope of dimension $d\leq 11$ is Ehrhart positive. We confirm that any order polytope of dimension $12$ or $13$ is Ehrhart positive. This solves an open problem proposed by Liu and Tsuchiya. Besides, we also verify that any $h^{*}$-polynomial of order polytope of dimension $d\leq 13$ is real-rooted.
\end{abstract}

\noindent
\begin{small}
 \emph{Mathematic subject classification}: Primary 05A15; Secondary 06A07, 68R05.
\end{small}

\noindent
\begin{small}
\end{small}

\section{Introduction}


A convex polytope is called \emph{integral} if any of its vertices has integer coordinates.
Let $\mathcal{P}$ represent a $d$-dimensional integral convex polytope in $\mathbb{R}^{n}$.
The function
$$\mathrm{ehr}(\mathcal{P},t)=|t\mathcal{P}\cap \mathbb{Z}^n|,\ \ \ \ \ t=1,2,\ldots,$$
counts the integer points within $t\mathcal{P}$, where $t\mathcal{P}=\{t\alpha : \alpha \in \mathcal{P}\}$ denotes the $t$-th dilation of $\mathcal{P}$.

Ehrhart \cite{Ehrhart62} proved that $\mathrm{ehr}(\mathcal{P},t)$ is a polynomial in $t$ of degree $d$.
This polynomial is referred to as the \emph{Ehrhart polynomial} of $\mathcal{P}$.
Moreover, the coefficient of $t^d$ in $\mathrm{ehr}(\mathcal{P},t)$ equals the (relative) volume of $\mathcal{P}$, the coefficient of $t^{d-1}$ is half of the boundary volume of $\mathcal{P}$, and the constant term is always $1$ (see \cite{BeckRobins}).
The remaining coefficients are complicated to describe \cite{McMullen77}.

An integral convex polytope $\mathcal{P}$ is said to be \emph{Ehrhart positive} (or have \emph{Ehrhart positivity}) if all the coefficients of $\mathrm{ehr}(\mathcal{P},t)$ are non-negative. In the case of dimension $3$, a well-known non-Ehrhart positive example is Reeve's tetrahedron \cite[Example 3.22]{BeckRobins}. For a comprehensive introduction to Ehrhart positivity, please refer to Liu's survey \cite{FuLiu19}. A recent development in Ehrhart positivity is presented in \cite{HibiHTY19}.

The generating function
$$\mathrm{Ehr}(\mathcal{P},x)=1+\sum_{t\geq 1}\mathrm{ehr}(\mathcal{P},t)x^t$$
is known as the \emph{Ehrhart series} of $\mathcal{P}$.
It is of the form
$$\mathrm{Ehr}(\mathcal{P},x)=\frac{h^{*}(x)}{(1-x)^{d+1}},$$
where $\dim\mathcal{P}=d$ and $h^{*}(x)$ is a polynomial in $x$ of degree $\deg h^{*}(x)\leq d$ (as detailed in \cite{BeckRobins}).
The $h^{*}(x)$ is usually referred to as the \emph{$h^{*}$-polynomial} of $\mathcal{P}$.

Stanley \cite{StanleyADM80} proved that the coefficients of the $h^{*}$-polynomial are nonnegative integers. Subsequently, plenty of work
has been done on various properties of the $h^{*}$-polynomial of an integral polytope. Let
$f(x)=a_mx^m+a_{m-1}x^{m-1}+\cdots+a_1x+a_0$ be a polynomial with nonnegative real coefficients.
The polynomial $f(x)$ is said to be \emph{unimodal} if its coefficients satisfy $a_0\leq \cdots \leq a_{i-1}\leq a_i\geq a_{i+1} \geq \cdots \geq a_m$ for some $0\leq i\leq m$. It is called \emph{log-concave} if $a_i^2\geq a_{i-1}a_{i+1}$ for $1\leq i\leq m-1$.
It is said to be \emph{real-rooted} if all roots of $f(x)$ are real numbers. If $f(x)$ is real-rooted, then it is log-concave;
If $f(x)$ is log-concave with $a_i>0$ for all $i$, then it is unimodal. See \cite[Section 5]{RP.StanleyAC} or \cite{StanleyLog-concave}.


This paper focuses on the families of order polytopes arising from finite partially ordered sets (or poset for short), which were initially introduced and studied by Stanley. See \cite{StanleyDCG} for detailed definitions.
\begin{dfn}[Stanley \cite{StanleyDCG}]\label{Definitorderp}
Given a poset $P$ on the set $[p]:=\{1,2,\ldots,p\}$, we associate a polytope $\mathcal{O}(P)$, called the \emph{order polytope} of $P$. It is defined to be the convex polytope consisting of those $(x_1,x_2,\ldots,x_p)\in \mathbb{R}^p$ such that
\begin{align*}
&x_i\leq x_j \ \ \ \text{if}\ \ \ i\prec_{P} j;
\\& 0\leq x_i\leq 1 \ \ \ \text{for}\ \ \ 1\leq i\leq p.
\end{align*}
\end{dfn}
The order polytope is an integral polytope of dimension $\dim \mathcal{O}(P)=p$ \cite{StanleyDCG}.

Regarding the Ehrhart positivity of order polytopes, Stanley provided an example of a non-Ehrhart-positive order polytope of dimension $21$ in \cite{StanleyMathOverfolw} (or \cite{Alexandersson}). Subsequently, Liu and Tsuchiya \cite{LiuTsuchiya19} proved that for any positive integer $d\geq 14$, there exists a non-Ehrhart-positive order polytope of dimension $d$. They also confirmed that any order polytope of dimension $d\leq 11$ is Ehrhart positive.
Naturally, they proposed the following open question.

\begin{prob}[Liu--Tsuchiya, \cite{LiuTsuchiya19}]\label{Open-Ehrhart}
For $d=12$ or $13$, does there exist an order polytope of dimension $d$ such that its Ehrhart polynomial has a negative coefficient?
\end{prob}

Our first contribution in this paper is to provide a negative solution to the above problem. That is, we confirm that any order polytope of dimension $d\leq 13$ is Ehrhart positive. 

Now the $h^{*}$-polynomial of order polytope $\mathcal{O}(P)$ is denoted by $h^{*}(x;P)$.
For any graded poset $P$, Reiner and Welker \cite{Reiner05} proved that the $h^{*}(x;P)$ is unimodal.
Br\"and\'en \cite{Branden05} showed that the $h^{*}(x;P)$ is not only unimodal but also symmetric.
Recently, a comprehensive research development of $h^{*}$-polynomials was presented by Ferroni and Higashitani in \cite{Ferroni23}.
They also raised the following open questions concerning order polytopes.

\begin{prob}{\em \cite[Section 6]{Ferroni23} or \cite[Section 7]{Stembridge09}}
For any graded poset $P$, prove that the $h^{*}$-polynomial $h^{*}(x;P)$ is real-rooted.
\end{prob}
There exists an order polytope whose $h^{*}$-polynomial is not real-rooted. Such an example was given by Stembridge in \cite{Stembridge09}. Similarly, Ferroni and Higashitani also raised the following open questions.

\begin{prob}{\em \cite[Section 6]{Ferroni23} or \cite[Section 7]{Stembridge09}}
For a general poset $P$, prove that the $h^{*}$-polynomial $h^{*}(x;P)$ is log-concave (or even unimodal).
\end{prob}

We verify that the $h^{*}$-polynomial of order polytopes $\mathcal{O}(P)$ of dimension $d\leq 13$ is real-rooted.
Consequently, these $h^{*}$-polynomials are log-concave and unimodal.

This paper is organized as follows.
In Section 2, we will briefly introduce how to compute the Ehrhart polynomial by Stembridge's \texttt{posets} package \cite{Stempackage}.
In Section 3, we first introduce how to generate posets. Then we study the Ehrhart positivity of order polytopes of dimensions $12$ and $13$. We also verify the real-rootedness of the $h^{*}(x;P)$ of $\#P\leq 13$.

\section{Ehrhart Polynomial and $h^{*}$-polynomial}

Let $P$ be a finite poset on $[p]$.
The $h^{*}$-polynomial of the order polytope $\mathcal{O}(P)$ is denoted by
$$h^{*}(x;P)=h_px^p+h_{p-1}x^{p-1}+\cdots+h_1x+h_0.$$
Equivalently, the Ehrhart polynomial of $\mathcal{O}(P)$ is given by 
$$\mathrm{ehr}(\mathcal{O}(P),t)=\sum_{i=0}^p h_i\binom{t+p-i}{p}.$$

Stembridge's \texttt{posets} package \cite{Stempackage} can be used to compute the Ehrhart polynomial of $\mathcal{O}(P)$, but not directly.
Denoted by $\mathbf{t}$ the $t$-element chain on the set $[t]$. The number of order-preserving maps from $P$ to $\mathbf{t}$ turns out to be a polynomial in $t$
of degree $p$. This polynomial is referred to as the \emph{order polynomial} of $P$, and denoted $\Omega_{P}(t)$.
Stanley \cite{StanleyDCG} showed that 
$$\mathrm{ehr}(\mathcal{O}(P),t)=\Omega_{P}(t+1).$$
The command \texttt{omega} in the \texttt{posets} package computes the order polynomial of $P$ directly, and hence 
can be used to compute the Ehrhart polynomial. 
The package \texttt{posets} is available at: https://dept.math.lsa.umich.edu/$\sim$jrs/.

There are other algorithms to compute the Ehrhart polynomial or the $h^*$ polynomial of a given order polytope. 
Stembridge's \texttt{posets} package seems to be the best for order polytopes of dimension no more than $13$ elements. 
Indeed, the command \texttt{omega} for computing the order polynomial has two options:
``\texttt{linear}" and ``\texttt{ideals}". The \texttt{linear} algorithm has minimal space requirements and a running time that is linear
in the number of linear extensions of $P$; The \texttt{ideals} algorithm is more appropriate for larger posets, and has worst-case time and space requirements that are quadratic in the number of order ideals of $P$. If neither algorithm is specified, the default is to use the \texttt{linear} option for posets with $\leq 6$ vertices, and the \texttt{ideals} option otherwise.

\section{Main Results}
To attack Open Problem \ref{Open-Ehrhart}, one needs only to generate posets $P$ of the desired dimension, compute the Ehrhart polynomial $\mathrm{ehr}(\mathcal{O}(P),t)$
and check its Ehrhart positivity. This turns out to take plenty of computation time. We first describe how to efficiently generate the posets, and then report our computation data.

\subsection{Generating Posets}
Brinkmann and Mckay \cite{Brinkmann02} determined the number of posets on up to $16$ points.
Table \ref{TabPoset} lists the number of posets for up to $13$ points.
The {notation $\text{poset}(p)$ refers} to the number of posets on the set $[p]$.

\begin{tiny}
\begin{table}[htbp]
    	\centering
    	\caption{The number of posets up to $13$ points}
    	\begin{tabular}{c||c|c|c|c|c|c|c|c}
    		\hline \hline
$p$ & 1 & 2 & 3 & 4 & 5 & 6 & 7 & 8   \\
    		\hline
$\text{poset}(p)$ & 1 & 2 & 5 & 16 & 63 & 318 & 2045 & 16999    \\
    		\hline
$p$ & 9 & 10&  11 &  12 & 13 & --- & --- & ---     \\
    		\hline
$\text{poset}(p)$ &  183231 & 2567284 &  46749427 & 1104891746 & 33823827452 & --- & --- & ---   \\
    		\hline
    	\end{tabular}\label{TabPoset}
\end{table}
\end{tiny}

\texttt{SageMath }\cite{SageMath} can generate all posets with a specified number of elements, but is not efficient.
We use the algorithm \texttt{nauty} in \cite{McKay13} for generating posets.
The C-package \texttt{nauty} \cite{Mckay2024} is available at: https://users.cecs.anu.edu.au/$\sim$bdm/nauty/.
One can use the procedure ``genposetg" to generate all posets with $p$ elements. The timings spent for $p=11, 12, 13$ are
approximately 2 seconds, 35 seconds, and 15 minutes, respectively. 
Note that this requires a large amount of storage space, as the number of posets is very large when $p=11,12,13$ (see Table \ref{TabPoset}).

\subsection{Ehrhart Positivity}
Stembridge's \texttt{posets} package contains all posets of at most $8$ elements. 
Stanley observed in \cite{StanleyMathOverfolw} that the order polytope of any poset with at most $8$ elements is Ehrhart positive. 

Liu and Tsuchiya \cite{LiuTsuchiya19} confirmed that the order polytope of any poset with at most $11$ elements is Ehrhart positive by using \texttt{SageMath} and Stembridge's \texttt{posets} package. It took Liu and Tsuchiya about three weeks to complete the computation for all posets of $11$ elements.
Through communication with Tsuchiya, we learned that most of the time was spent for generating the posets.

We use C-package \texttt{nauty} to generate the posets of $12$ elements and $13$ elements. Then we transform the output posets into their Maple format. Finally, we use the Maple-package \texttt{posets} to verify the Ehrhart positivity of order polytopes of dimension $12$ and $13$.

Based on \cite[Theorem 1.7]{LiuTsuchiya19}, we summarize our findings as follows.
\begin{thm}\label{Ehrhart13-14}
Any order polytope of dimension $d\leq 13$ is Ehrhart positive. For any positive integer $d\geq 14$, there exists a non-Ehrhart positive order polytope of dimension $d$.
\end{thm}

\subsection{Real-Rootedness}

A poset is defined as \emph{narrow} if its vertices can be partitioned into two chains. According to Dilworth's Theorem, a poset is narrow if and only if it contains no antichain of $3$ elements.
Stembridge \cite{Stembridge09} provided an example of a narrow poset with $17$ vertices whose $h^{*}$-polynomial has nonreal zeros.
This disproves Neggers' conjecture \cite{Neggers} that the $h^{*}$-polynomial of $\mathcal{O}(P)$ have all real zeros.
Stembridge also proved that the $h^{*}$-polynomial of the order polytope of a narrow poset with at most $16$ vertices is real rooted.
In \cite{Stembridge09}, Stembridge verified that if $P$ is a poset on at most $10$ vertices, then its $h^*$ polynomial $h^{*}(x;P)$ is real rooted.

We employ Sturm's Theorem (for example, see \cite[Page 419]{Knuth81} or \cite{Sturm29}) to calculate the number of real roots of a polynomial.

Let $f(x)=a_mx^m+a_{m-1}x^{m-1}+\cdots+a_1x+a_0$ be a polynomial of degree $m$ with real coefficients.
The \emph{Sturm sequence} of the polynomial $f(x)$ is the sequence of polynomials $(f_0(x), f_1(x),\ldots,)$ recursively defined by
$$f_0(x)=f(x),\ \ \ f_1(x)=f^{\prime}(x),\ \ \ f_{i+1}(x)=-\texttt{rem}(f_{i-1}(x),f_i(x)),$$
for $i\geq 1$, where $f^{\prime}(x)$ is the derivative of $f(x)$, and $\texttt{rem}(f_{i-1}(x),f_i(x))$ is the remainder of
$f_{i-1}(x)$ when divided by $f_i(x)$. The length of the Sturm sequence is at most the degree of $f(x)$.

The number of sign variations at $c\in\mathbb{R}$ of the Sturm sequence of $f(x)$ is the number of sign changes (ignoring zeros) in the sequence of real numbers $f_0(c), f_1(c), f_2(c),\ldots$. The number of sign variations is denoted by $V(c)$.

\begin{thm}{\em \cite{Sturm29}}\label{Sturm-Theorem}
Let $f(x)$ be a univariate polynomial with real coefficients. Let $a,b\in \mathbb{R}$ with $a<b$, and $f(a),f(b)\neq 0$.
If $f(x)$ has distinct roots (i.e., $\gcd(f(x), f^{\prime}(x))=1$), then the number of real roots of $f(x)$ in $(a,b]$ is $V(a)-V(b)$.
\end{thm}
Letting $a\to -\infty$ and $b\to \infty$ give what we want. Clearly, the sign of $f(\infty)$ is the sign of its leading coefficient $a_m$,
and the sign of $f(-\infty)$ is the sign of $a_m$ if $m$ is even and the sign of $-a_m$ if $m$ is odd.

We apply Theorem \ref{Sturm-Theorem} to verify the real rootedness of the $h^{*}$-polynomial of the order polytope $\mathcal{O}(P)$ for dimension $d\leq 13$. This result was validated using \texttt{Maple}. We summarize the results as follows.

\begin{thm}
The $h^{*}$-polynomial of any order polytope of dimension $d\leq 13$ is real-rooted,
and is hence log-concave and unimodal.
\end{thm}

\subsection{Timing}

In this paper, we employ \texttt{Maple 2021} \cite{Maple} to confirm the above results.
That is, the Ehrhart positivity and the real-rootedness of $h^{*}(x;P)$ of order polytopes of dimension $\leq 13$.
Utilizing a CPU with a frequency of $2.7$ GHz, the verification of Ehrhart positivity and real-rootedness of $h^{*}(x;P)$ for dimensions $12$ and $13$ required $3160$ and $168669$ CPU hours, respectively. To expedite these computations, we used a parallel processing approach with $192$ CPUs. The actual running time
for the verifications is approximately $37$ days.

For the reader's convenience, we provide all the details of our calculations in the supplementary electronic material.
The \texttt{Maple} program can be freely obtained at the first author'swebsite:\\ https://github.com/TygerLiu/TygerLiu.github.io/tree/main/Procedure/orderpolytope.




\noindent
{\small \textbf{Acknowledgements:}
We are grateful to Prof. Fu Liu and Prof. Akiyoshi Tsuchiya for many useful discussions.
We also are indebted to Prof. Brendan D. Mckay for helpful advice about the package \textbf{nauty}.
The authors would like to thank the anonymous referee for valuable suggestions for improving the presentation.
This work is partially supported by the National Natural Science Foundation of China [12071311].

\end{document}